\documentclass[10pt]{article}
\usepackage[final]{epsfig}
\usepackage{graphics}
\usepackage{amsmath}
\usepackage{amsfonts}
\usepackage{latexsym}
\usepackage{amssymb}
\usepackage{graphicx}
\usepackage[all,poly,knot,dvips]{xy}
\usepackage{amscd}
\usepackage{amsthm}
\usepackage{url}

\newtheorem{theorem}{Theorem}[section]

 \newtheorem{proposition}[theorem]{Proposition}
 \newtheorem{corollary}[theorem]{Corollary}

\newtheorem{remark}[theorem]{Remark}

\begin{document}
\newcommand{\eps}{{\varepsilon}}
\newcommand{\g}{{\gamma}}
\newcommand{\G}{{\Gamma}}
\newcommand{\proofend}{$\Box$\bigskip}
\newcommand{\C}{{\mathbf C}}
\newcommand{\Q}{{\mathbf Q}}
\newcommand{\R}{{\mathbf R}}
\newcommand{\Z}{{\mathbf Z}}
\newcommand{\RP}{{\mathbf {RP}}}
\newcommand{\CP}{{\mathbf {CP}}}
\newcommand{\Tr}{{\rm Tr\ }}
\def\proof{\paragraph{Proof.}}
\def\one{\mathbf 1}
\def\bop{{\mathbf p}}
\def\boq{{\mathbf q}}
\def\boxx{{\mathbf x}}
\def\boy{{\mathbf y}}
\def\barT{{\overline{T}}}
\def\barD{{\overline{D}}}
\def\caq{{\cal Q}}
\def\caP{{\cal P}}
\def\can{{\cal N}}
\def\EE{{\cal E}}
\def\FF{{\cal F}}
\newcommand{\PSL}{{\mathrm{PSL}}}
\newcommand{\PGL}{{\mathrm{PGL}}}
\newcommand{\SL}{{\mathrm{SL}}}

\title{Higher pentagram maps, weighted directed networks, and cluster dynamics}
\author{Michael Gekhtman,\thanks{
Department of Mathematics, University of Notre Dame, Notre Dame, IN 46556, USA;
e-mail: \tt{mgekhtma@nd.edu}
}
\ Michael Shapiro,\thanks{
Department of Mathematics, Michigan State University, East Lansing, MI 48823, USA;
e-mail: \tt{mshapiro@math.msu.edu}
}
\  Serge Tabachnikov,\thanks{
Department of Mathematics,
Pennsylvania State University, University Park, PA 16802, USA;
e-mail: \tt{tabachni@math.psu.edu}
}\\
\ and Alek Vainshtein\thanks{
Department of Mathematics and Department of Computer Science, University of Haifa, Haifa, Mount 
Carmel 31905, Israel;
e-mail: \tt{alek@cs.haifa.ac.il}
}
\\
}
\maketitle

\begin{abstract}
The pentagram map was extensively studied in a series of papers by by V. Ovsienko, R. Schwartz and S. Tabachnikov. It was recently interpreted by  M. Glick  as a sequence
of cluster transformations associated with a special quiver. Using compatible Poisson structures in cluster algebras and Poisson geometry of directed networks
on surfaces, we generalize Glick's construction 
to include the pentagram map into a  family of geometrically meaningful  discrete integrable maps.
\end{abstract}

\section{Introduction} \label{intro}

\begin{figure}[hbtp]
\centering
\includegraphics[height=1.5in]{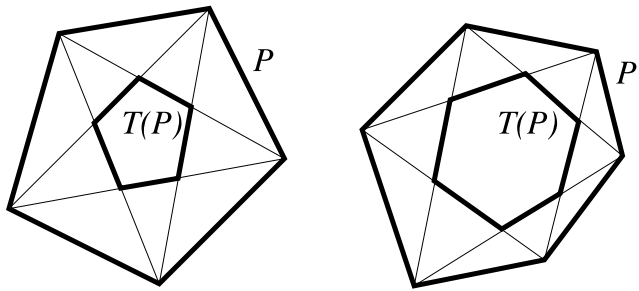}
\caption{Pentagram map}
\label{Penta}
\end{figure}

The pentagram map was introduced by R. Schwartz about 20 years ago \cite{Sch1}. Recently, it has attracted a considerable attention: see \cite{Gl,MB,MOT,OST1,OST2,OST3,Sch2,Sch3,ST1,ST2,So} 
for various aspects of the pentagram map and related topics. 
On plane polygons, the pentagram map (depicted in Fig.~\ref{Penta}) acts  by drawing the diagonals that connect second-nearest vertices of a poygon $P$ and forming a new polygon $T(P)$ whose vertices are their consecutive intersection points. The pentagram map commutes with projective transformations, so it acts on the projective equivalence classes of polygons in the projective plane.

The pentagram map can be extended to a larger class of {\it twisted polygons}. A twisted $n$-gon is a sequence of points $V_i \in RP^2$ in general position such that $V_{i+n}=M(V_i)$ for all $i \in \Z$ and some fixed element $M\in\PGL(3,\R)$, called the monodromy. 
A polygon is closed if the monodromy is the identity. 
Denote by ${\cal P}_n$ the 
space of projective equivalence classes of twisted $n$-gons; this is a variety of dimension $2n$. Denote by  $T:{\cal P}_n\to {\cal P}_n$  the pentagram map. (The $i$th vertex of the image is the intersection of diagonals $(V_i,V_{i+2})$ and $(V_{i+1},V_{i+3})$.)

One has a coordinate system $X_1,Y_1,\dots,X_n,Y_n$ in ${\cal P}_n$ where $X_i,Y_i$ are the so-called corner invariants associated with $i$th vertex, discrete versions of  projective curvature, see \cite{Sch3,OST1,OST2}. In these coordinates, the pentagram map is a rational transformation
\begin{equation} \label{pentaformula}
X^*_i=X_i\,\frac{1-X_{i-1}\,Y_{i-1}}{1-X_{i+1}\,Y_{i+1}},
\qquad 
Y^*_i=Y_{i+1}\,\frac{1-X_{i+2}\,Y_{i+2}}{1-X_{i}\,Y_{i}}
\end{equation}
(the indices are taken mod $n$). Schwartz \cite{Sch3} 
observed that the pentagram map commutes with the $\R^*$ action on ${\cal P}_n$ given by
$(X_i,Y_i)\mapsto (tX_i, t^{-1}Y_i)$.

The main feature of the pentagram map is that it is a discrete completely integrable system (its continuous $n\to \infty$ limit is the Boussinesq equation, one of the best known completely integrable PDEs). 
Specifically, the pentagram map has $2[n/2]+2$ integrals constructed in \cite{Sch3}; these integrals are polynomial in the $X,Y$-coordinates and they are algebraically independent. The space ${\cal P}_n$ has a $T$-invariant Poisson structure, introduced in \cite{OST1,OST2}. The corank of this Poisson structure equals 2 or 4, according as $n$ is odd or even, and the integrals are in involution. This provides Liouville integrability of the pentagram map on the space of twisted polygons. 

Complete integrability on the smaller space ${\cal C}_n$ of closed $n$-gons is proved in \cite{OST3}. F. Soloviev \cite{So} established algebraic-geometric integrability of the pentagram map by constructing its Lax (zero curvature) representation; his approach works both for twisted and closed polygons. 

It was shown in \cite{Sch3} that the pentagram map was intimately related to the so-called octahedral recurrence. It was conjectured in \cite{OST1,OST2} that
the pentagram map was related to cluster
transformations. This relation was discovered and explored by Glick  
[12]
who proved that the pentagram map, acting on the quotient space ${\cal P}_n/\R^*$,
is described by coefficient dynamics \cite{FZ4} (also known as 
$\tau$-transformations, see Chapter 4 in \cite{GSV3}) for a certain cluster
structure. 

In this research announcement and the forthcoming detailed paper, we extend and generalize Glick's work by including the pentagram map into a family of discrete completely integrable systems. 
Our main tool is Poisson geometry of weighted  directed networks on surfaces.
The ingredients necessary for complete integrability -- invariant Poisson brackets, integrals of motion in involution, Lax representation -- are recovered from combinatorics of the networks.
Postnikov \cite{Po} introduced such networks 
in the case of a disk and investigated their transformations and their relation  to cluster
transformations; most of his results are local, and
hence remain valid for networks on any surface. Poisson  
properties of
weighted directed networks in a disk and their relation to r-matrix
structures on $GL_n$ are studied in \cite{GSV2}. In \cite{GSV4} these results were
further extended to networks in an annulus and r-matrix Poisson  
structures on
matrix-valued rational functions. Applications of these techniques  
to the
study of integrable systems can be found in \cite{GSV5}. A detailed  
presentation
of the theory of weighted directed networks from a cluster algebra
perspective can be found in  Chapters 8--10 of \cite{GSV3}.

Our integrable systems depend on one discrete parameter $k\ge2$. The  
case $k=3$ corresponds to the pentagram map. For $k\ge4$,
we give our integrable systems a geometric interpretation as pentagram-like maps involving deeper diagonals. If $k=2$ and the ground field is $\C$, we give a geometric interpretation in terms of circle patterns \cite{Sc,BH}.

While working on this manuscript, we were informed by V.~Fock that in his current project with A. Marshakov, closely related integrable systems are studied using the approach of \cite{GoK}.

\section{Generalized Glick's quivers and the\\ $(\bop,\boq)$-dynamics} \label{pq}

\begin{figure}[hbtp]
\centering
\includegraphics[height=1in]{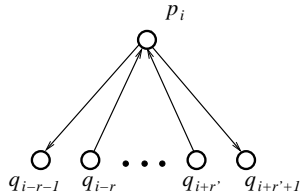}
\caption{The quiver ${\cal Q}_{k,n}$}
\label{quiver}
\end{figure}

For any integer $n\ge 2$ and any integer $k$, $2\le k\le n$,  
consider the quiver (an oriented multigraph without loops) ${\cal Q}_{k,n}$ 
defined as follows: $\caq_{k,n}$ is a
bipartite graph on $2n$ vertices labeled $p_1,\dots, p_n$ and $q_1,\dots,q_n$ (the labeling is cyclic, so that $n+1$ is the same as $1$). The graph is 
invariant under the shift $i\mapsto i+1$.
Each vertex has two incoming and two outgoing edges. The number $k$ is the ``span" of the quiver, the distance between two outgoing edges from a $p$-vertex, see Fig.~\ref{quiver} where 
$r=[k/2]-1$ and $r+r'=k-2$ (in other words, $r'=r$ for $k$ even and $r'=r+1$ for $k$ odd). For $k=3$, we have Glick's  quiver \cite{Gl}.

We consider the cluster structure associated with the quiver ${\cal Q}_{k,n}$. Choose variables $\bop=(p_1,\dots,p_n)$ and
$\boq=(q_1,\dots,q_n)$ as $\tau$-coordinates (see \cite{GSV3}, Chapter 4),  
and consider cluster transformations corresponding to the quiver mutations at the $p$-vertices. These transformations commute, and we perform them simultaneously. This leads to the transformation (the new variables are marked by asterisque)
\begin{equation} \label{exch}
p_i^*=\frac{1}{p_i},\quad q_i^*=q_i\frac{(1+p_{i-r'-1})(1+p_{i+r+1})p_{i-r'}p_{i+r}}{(1+p_{i-r'})(1+p_{i+r})},
\end{equation} 
(see Lemma 4.4 of \cite{GSV3} for the exchange relations for $\tau$-coordinates). 
The resulting  quiver is identical to $Q_{k,n}$ with the letters $p$ and $q$ interchanged. Thus we compose transformation (\ref{exch}) with the involution $\{p_i \leftrightarrow q_i\}_{i=1}^n$ for $k$ even,
or with the transformation $\{p_i\mapsto q_{i+1}, q_i\mapsto p_i\}_{i=1}^n$ for $k$ odd, and arrive at the transformation that we denote by $\overline{T}_k$:
\begin{equation} \label{mappq}
\begin{split} 
q_i^*=\frac{1}{p_i},\quad p_i^*=q_i\frac{(1+p_{i-r-1})(1+p_{i+r+1})p_{i-r}p_{i+r}}{(1+p_{i-r})(1+p_{i+r})},\quad k\ {\rm even},\\
q_i^*=\frac{1}{p_{i-1}},\quad p_i^*=q_i\frac{(1+p_{i-r-2})(1+p_{i+r+1})p_{i-r-1}p_{i+r}}{(1+p_{i-r-1})(1+p_{i+r})},\quad k\ {\rm odd}.
\end{split} 
\end{equation} 
The difference in the formulas is due to the asymmetry between left and right in the enumeration of vertices 
in Fig.~\ref{quiver} for odd $k$, when $r'\ne r$.

Let us equip the $(\bop,\boq)$-space with a Poisson structure compatible with the cluster structure,
see~\cite{GSV3}. Denote by $A=(a_{ij})$ the $2n\times 2n$ skew-adjacency matrix of ${\cal Q}_{k,n}$, assuming that the first $n$ rows and columns correspond to $p$-vertices. Then we put $\{v_i,v_j\}=a_{ij}v_iv_j$, where $v_i=p_i$ for $1\le i\le n$ and $v_i=q_{i-n}$ for $n+1\le i\le 2n$.

\begin{theorem}\label{pqsyst}
{\rm (i)} The above Poisson structure is invariant under the map $\overline{T}_k$. 

{\rm (ii)} The function $\prod p_iq_i$ is an integral of the map $\overline{T}_k$. Besides, it is
Casimir, and hence the Poisson structure and the map descend to the hypersurfaces $\prod p_iq_i=\rm{const}$.
\end{theorem}
 
We denote by $\overline{T}^c_k$ the restriction of $\overline{T}_k$ to the hypersurface $\prod p_iq_i=c$.
In what follows, we shall be concerned only with $\overline{T}^1_k$,  which we shorthand to $\overline{T}_k$.
Note that $\overline{T}^1_3$ is the  pentagram map on ${\cal P}_n/\R^*$ considered by Glick \cite{Gl}.

Let us consider an auxiliary transformation $\barD_k$ given by $\{p_i\mapsto 1/q_i, q_i\mapsto 1/p_i\}_{i=1}^n$ for $k$ even
and $\{p_i\mapsto 1/q_{i+1}, q_i\mapsto 1/p_i\}_{i=1}^n$ for $k$ odd. Then $\barT_k$ and its inverse are related via
\begin{equation}\label{barDTD}
\barD_k\circ\barT_k\circ\barD_k=\barT_k^{-1}.
\end{equation}

\begin{figure}[hbtp]
\centering
\includegraphics[height=0.9in]{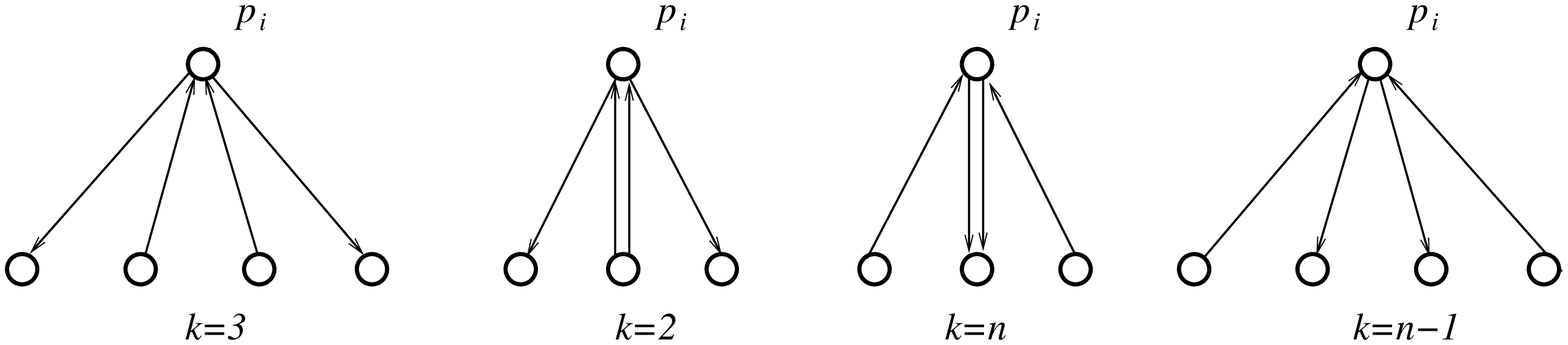}
\caption{The quivers ${\cal Q}_{k,n}$ for various values of $k$}
\label{pandq}
\end{figure}

\begin{remark}
{\rm Along with the $p$-dynamics when the mutations are performed at the $p$-vertices of the quiver ${\cal Q}_{k,n}$, one may consider the respective $q$-dynamics. We note that the $q$-dynamics is essentially the same as $p$-dynamics. Namely, the $q$-dynamics for a given value of $k$ corresponds to the $p$-dynamics for $k'=n+2-k$. This is illustrated in Fig.~\ref{pandq}.} 
\end{remark} 

\section{Weighted directed networks and the \\
$(\boxx,\boy)$-dynamics} \label{xysect}
 
{\bf Weighted directed networks on surfaces.}
We start with a very brief description of the theory of weighted directed networks on surfaces with a boundary, adapted for our purposes; see \cite{Po,GSV3} for details. In this note, we will only need to consider acyclic graphs on a cylinder (equivalently, annulus) $\mathcal{C}$ that we position horizontally with two boundary circles, one on the left and one on the right.

 Let $G$ be a directed acyclic graph with the vertex set $V$ and the edge set $E$ embedded in $\mathcal{C}$. 
$G$ has $2n$ 
{\it boundary vertices\/}, each of degree one: $n$ {\it sources\/} on the left boundary circles and $n$ {\it sinks\/} on the right
boundary circle.
Both sources and sinks are numbered clockwise
as seen from the axis of a cylinder behind the left boundary circle.
All the internal vertices of $G$ have degree~$3$ and are of two types: either they have exactly one
incoming edge (white vertices), or exactly one outgoing edge (black vertices).  
To each edge $e\in E$ we assign the {\it edge weight\/} $w_{e}\in\R\setminus 0$. 
A {\it perfect network\/} $\can$ is obtained from $G$  
by adding an oriented  curve $\rho$ without self-intersections (called a {\it cut\/}) that joins  the left
and the right boundary circles and does not contain vertices of $G$. The points of the {\it space of edge weights\/} $\EE_\can$
can be considered as copies of $\can$  with edges weighted by nonzero real numbers.

Assign an independent variable $\lambda$ to the cut $\rho$. The weight of a directed path $P$ between a source and a sink 
is defined as the product of the weights of all edges along the path times  $\lambda$ raised 
into the power equal to the intersection index of $\rho$ and $P$ (we assume that all intersection points are transversal, in which case the intersection index is the number of intersection points counted with signs).
The {\it boundary measurement\/} between $i$th source  and  $j$th sink is then defined as the sum of path weights over all 
(not necessary simple) paths between them. A boundary measurement is rational
 in the weights of edges and $\lambda$, see \cite{GSV4}.

Boundary measurements are organized in the {\it boundary measurement matrix}, thus giving rise to the {\it boundary measurement map\/}
from $\EE_\can$ to the space of $n\times n$ rational matrix 
functions. The gauge group acts on $\EE_\can$ as follows: for any internal vertex $v$ of $\can$ and any Laurent monomial $L$ 
in the weights $w_e$ of $\can$, the weights of all edges leaving $v$ are multiplied by $L$, and the weights of all edges entering $v$ are
 multiplied by $L^{-1}$. Clearly, the weights of paths between boundary vertices, and hence boundary measurements, 
 are preserved under this action.
Therefore, the boundary measurement map can be factorized through the space $\FF_\can$ defined as the 
quotient of $\EE_\can$ by the action of the gauge group. 

It is explained in \cite{GSV4} that $\FF_\can$ can be parametrized as follows.
The graph $G$ divides $\mathcal C$ into a finite number of connected components called
\emph{faces}. The boundary of each face consists of edges of $G$ and, possibly, of several arcs of 
$\partial \mathcal C$. 
A face is called {\it bounded\/} if its boundary contains only edges of $G$ and {\it unbounded\/} otherwise. 
Given a face $f$, we define its {\it face weight\/} 
$y_f=\prod_{e\in\partial f}w_e^{\gamma_e}$,
where $\gamma_e=1$ if the direction of $e$ is compatible with the counterclockwise orientation of the 
boundary $\partial f$ and $\gamma_e=-1$ otherwise. 
Face weights are invariant under the gauge group action. 
Then $\FF_\can$ is parametrized by the collection of all face weights and a weight of a fixed path in $G$ joining two boundary circles.

\begin{figure}[hbtp]
\centering
\includegraphics[height=2.3in]{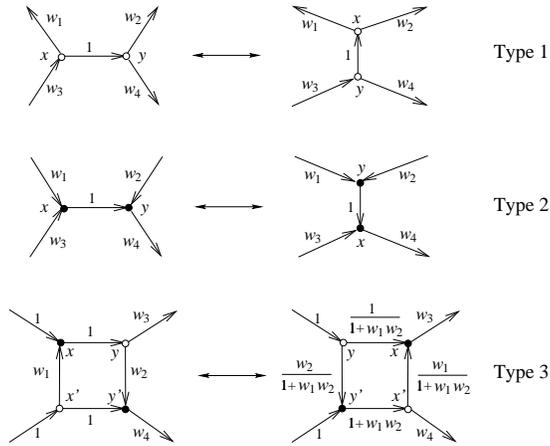}
\caption{Postnikov transformations }
\label{moves}
\end{figure}

Below we will frequently use elementary 
transformations of weighted networks that do not change
the boundary measurement matrix. They were introduced by Postnikov in \cite{Po} and are presented in Figure \ref{moves}.

As was shown in \cite{GSV2, GSV4}, the space of edge weights can be made into a Poisson manifold 
by considering Poisson brackets 
that behave nicely with respect to a natural operation of concatenation of networks. Such Poisson brackets on $\EE_\can$ form a 6-parameter
family, which is pushed forward to a 2-parameter family of Poisson brackets on $\FF_\can$. Here we will need a  specific member of the latter family. The corresponding Poisson structure, called {\em standard}, is described in terms of
the {\it directed dual network\/} $\can^*$ defined as follows. Vertices of $\can^*$ are the faces of $\can$. 
Edges of $\can^*$ correspond to the edges of $\can$ that connect either two internal vertices of different colors, 
or an internal vertex with a boundary vertex; note that there might be several edges between the same pair of 
vertices in $\can^*$. An edge $e^*$ in $\can^*$ corresponding to $e$ in $\can$ is directed in such a way that the white endpoint of $e$ (if it exists) lies to the left of $e^*$ and 
the black endpoint of $e$ (if it exists) lies to the right of $e$. 
The weight $w^*(e^*)$ equals $1$ if both endpoints of $e$ are internal vertices, and $1/2$ if one of the 
endpoints of $e$ is a boundary vertex.  Then
the restriction of the standard Poisson bracket on $\FF_\can$ to the space of face weights is given by
\begin{equation}
\label{facebracket}
\{y_f,y_{f'}\}=\left(\sum_{e^*: f\to f'} w^*(e^*)-
\sum_{e^*: f'\to f} w^*(e^*)\right)y_fy_{f'}.
\end{equation}

Any network $\can$ of the kind described above gives rise to a network $\bar\can$ on a torus. To do this, one identifies boundary circles in such a way that every sink is glued to a source with the same label. The resulting two-valent vertices are then erased, so that every pair of glued edges becomes
a new edge with the weight equal to the product of two edge-weights involved. Similarly, $n$ pairs
of unbounded faces are glued together into $n$ new faces, whose face-weights are products
of pairs of face-weights involved. We will view two networks on a torus as equivalent if their 
underlying graphs differ only by orientation of edges, but have the same vertex coloring and the same face
weights. The parameter space we associate with $\bar\can$ consists of face weights and the weight $x$ of a single fixed
directed cycle homological to the closed curve on the torus obtained by identifying endpoints of the cut.
The standard Poisson bracket induces a Poisson bracket on face-weights of the new network, which is again given by~\eqref{facebracket} with the dual graph $\can^*$ replaced by $\bar\can^*$ defined by the same rules. The bracket between $x$ and face-weights is given by $\{x, y_f\} = c_f x y_f$, where $c_f$ is a constant whose exact value will not be important to us.

\begin{figure}[hbtp]
\centering
\includegraphics[height=0.75in]{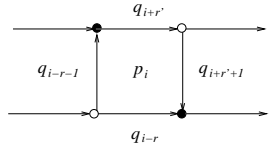}
\caption{Local structure of the network $\bar{\cal N}_{k,n}$  on the torus}
\label{names}
\end{figure}

{\bf The $(\boxx,\boy)$-dynamics.} 
Observe first of all that ${\cal Q}_{k,n}$ can be embedded in a unique way into a torus. Following \cite{GSV2}, we consider
a network ${\cal N}_{k,n}$ on the cylinder such that ${\cal Q}_{k,n}$ is the directed dual of the corresponding network $\bar\can_{k,n}$ on the torus. 
Applying Postnikov's transformations of types 1 and 2, 
one gets a network whose faces are quadrilaterals ($p$-faces) and octagons ($q$-faces). Locally, the network $\bar{\cal N}_{k,n}$ is shown in Fig.~\ref{names}.
Globally, $\bar\can_{k,n}$ consists of $n$ such pieces glued together in such a
way that the lower right edge of the $i$-th piece is identified with the
upper left edge of the $(i+1)$-th piece, and the upper right edge of the
$i$-th piece is identified with the lower left edge of the $(i+k-1)$-st
piece (all indices are mod $n$).

\begin{figure}[hbtp]
\centering
\includegraphics[height=1in]{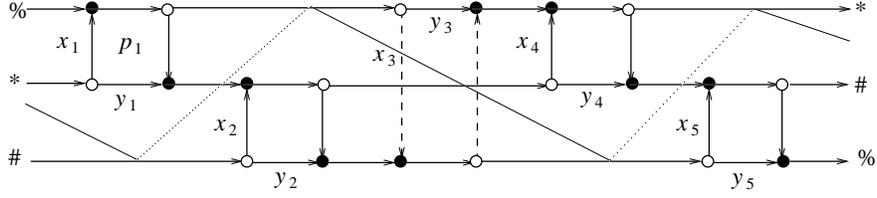}
\caption{The network $\bar{\cal N}_{3,5}$ on the torus}
\label{network}
\end{figure}

The network $\bar{\cal N}_{3,5}$  is shown in Fig.~\ref{network}. The figure depicts a torus, represented as a flattened two-sided  cylinder (the dashed lines are on the ``invisible" side); the edges marked by the same symbol are glued together accordingly. 
The cut is shown by the thin line.

\begin{figure}[hbtp]
\centering
\includegraphics[height=0.9in]{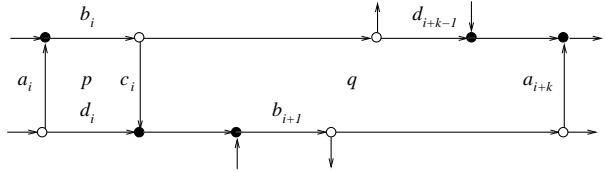}
\caption{Edge weights prior to the gauge group action}
\label{gauge}
\end{figure}

Assume that the edge weights around the face $p_i$ are $a_i$, $b_i$, $c_i$, and $d_i$; without loss of generality, we may assume that all other weights are equal~$1$, see Fig.~\ref{gauge}.

Applying the gauge group action, we can set to~$1$ the weights of the upper and the right edges of each quadrilateral face, while keeping weights of all edges with both endpoints of the same color equal to $1$.
For the face $p_i$, denote by $x_i$ the weight of the left edge and by $y_i$, the weight of the lower edge after the gauge group action (see Fig.~\ref{network}). Put
$\boxx=(x_1,\dots,x_n)$, $\boy=(y_1,\dots,y_n)$.

\begin{proposition}\label{weights}
{\rm (i)} The weights $(\boxx,\boy)$ are given by
$$ 
y_i=\frac{d_i}{b_ib_{i-1}\dots b_{i-k+2}c_ic_{i-1}\dots c_{i-k+1}},\quad x_i=\frac{a_i}{b_{i-1}\dots b_{i-k+2}c_{i-1}\dots c_{i-k+1}}.
$$

{\rm (ii)} The relation between $(\bop,\boq)$ and $(\boxx,\boy)$ is as follows:
$$
p_i=\frac{y_i}{x_i},\quad q_i=\frac{x_{i+r+1}}{y_{i+r}};\qquad x_i=x_1 \prod_{j=1}^{i-1}p_jq_{j-r},\quad y_i=x_i p_i.
$$
\end{proposition}

Note that the projection $\pi:(\boxx,\boy)\mapsto (\bop,\boq)$ has a 1-dimensional fiber.

The map $\overline{T}_k$ can be described via equivalent transformations of the network $\bar\can_{k,n}$.
The transformations include Postnikov's moves of types 1, 2, and 3, and the gauge group action.
 We describe the sequence of these transformations below.

\begin{figure}[hbtp]
\centering
\includegraphics[height=1in]{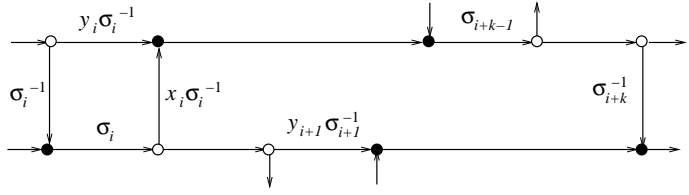}
\caption{Type 3 Postnikov's move for $\bar\can_{k,n}$}
\label{aftermut}
\end{figure}

First, we apply Postnikov's type 3 move at each $p$-face (this corresponds to cluster $\tau$-transformations at $p$-vertices
of $\caq_{k,n}$ given by~\eqref{exch}). Locally, the result is shown in Fig.~\ref{aftermut} where $\sigma_i=x_i+y_i$.

Next, we apply type 1 and type 2 Postnikov's moves at each white-white and black-black edge, respectively.
 In particular, we move vertical arrows interchanging the right-most and the left-most position on the network in 
 Fig.~\ref{network} using the fact that it is drawn on the torus. These moves interchange the quadrilateral and octagonal faces of the graph thereby swapping  the variables $p$ and~$q$.

It remains to use gauge transformations to achieve the weights  as in Fig.~\ref{network}. In our situation, weights  $a_i$, $b_i$, $c_i$, $d_i$ are as follows, see Fig.~\ref{aftermut}:
$$
a_i=\frac{x_i}{\sigma_i},\quad b_i=\sigma_{i+k-1},\quad c_i=\frac{1}{\sigma_{i+k}},\quad d_i=\frac{y_{i+1}}{\sigma_{i+1}}.
$$
Using Proposition~\ref{weights}(i), we obtain the new values of $(\boxx,\boy)$; 
we also shift the indices to conform with Fig.~\ref{names}.
This yields the map $T_k$, the main character of this note, described in the following proposition.

\begin{proposition} {\rm (i)} The map $T_k$ is given by
\begin{equation} \label{mapxy}
x_i^*=x_{i-r'-1} \frac{x_{i+r}+y_{i+r}}{x_{i-r'-1}+y_{i-r'-1}},\quad y_i^*=y_{i-r'} \frac{x_{i+r+1}+y_{i+r+1}}{x_{i-r'}+y_{i-r'}},
\end{equation}

{\rm (ii)} The maps $T_k$ and $\overline{T}_k$ are conjugated via $\pi$: $\pi\circ T_k=\barT_k\circ \pi$.
\end{proposition}

 Note that the map $T_k$ commutes with the scaling action of the group $\R^*$: $(\boxx,\boy)\mapsto (t\boxx,t\boy)$, and that the orbits of this action are the fibers of the projection $\pi$.
  
The map $T_3$ coincides with the pentagram map. Indeed, for $k=3$,~\eqref{mapxy} gives 
\begin{equation} \label{k=3}
x_i^*=x_{i-2} \frac{x_{i}+y_{i}}{x_{i-2}+y_{i-2}},\quad y_i^*=y_{i-1} \frac{x_{i+1}+y_{i+1}}{x_{i-1}+y_{i-1}}.
\end{equation}
Change the variables as follows:
$
x_i\mapsto Y_i,\ y_i\mapsto -Y_iX_{i+1}Y_{i+1}.
$
In the new variables, the map~\eqref{k=3} is rewritten as
$$
X_i^*=X_{i-1} \frac{1-X_{i-2}Y_{i-2}}{1-X_iY_i},\quad Y_i^*=Y_i \frac{1-X_{i+1}Y_{i+1}}{1-X_{i-1}Y_{i-1}},
$$
which differs from~\eqref{pentaformula} only by the cyclic shift of the index $i\mapsto i+1$. Note that the maps $T_k$, and in particular the pentagram map, commute with this shift.

The map $T_2$ is a periodic version of the discretization
of the relativistic Toda lattice suggested in \cite{Su}.
It belongs to a family of Darboux-B\"acklund transformations of integrable lattices of Toda type, that were put into a cluster algebras framework in \cite{GSV5}.

Let us define an auxiliary map $D_k$ by
\begin{equation}\label{DTD}
x_i^*=\frac{y_{i-r}y_{i-r+1}\cdots y_{i+r'-1}}{x_{i-r}x_{i-r+1}\cdots x_{i+r'}},\qquad
y_i^*=\frac{y_{i-r}y_{i-r+1}\cdots y_{i+r'}}{x_{i-r}x_{i-r+1}\cdots x_{i+r'+1}}.
\end{equation}

\begin{proposition}
{\rm (i)} The maps $T_k$ and $T_k^{-1}$ are related by
$$
D_k\circ T_k\circ D_k=T_k^{-1}.
$$

{\rm (ii)} The maps $D_k$ and $\barD_k$ are conjugated via $\pi$: $\pi\circ D_k=\barD_k\circ \pi$.

\end{proposition}

\section{Poisson structure and complete integrability} \label{integrability}

The main result of this note is {\em complete integrability} of                    
transformations $T_k$, i.e., the existence of a $T_k$-invariant Poisson  
bracket and of 
a maximal family of integrals in involution. The key ingredient of the proof is the result  
obtained
in \cite{GSV4} on Poisson properties of the boundary measurement map  
defined in Section~3.1. First, we recall the definition of an R-matrix  
(Sklyanin)
bracket, which plays a crucial role in the modern theory of integrable  
systems \cite{OPRS, FT}. The bracket is defined on the space of $n\times n$ rational matrix
functions $M(\lambda)=(m_{ij}(\lambda))_{i,j=1}^n$ and is given by  
the formula
\begin{equation}
\label{sklya}
\left \{M(\lambda){\stackrel{\textstyle{\small{\otimes}}}{,}} M(\mu)\right\}  
= \left [ R(\lambda,\mu), M(\lambda)\otimes M(\mu) \right ],
\end{equation}
where the left-hand is  understood as $\left \{M(\lambda) 
{\stackrel{\textstyle{\otimes}}{,}} M(\mu)\right\}_{ii'}^{jj'}=\{m_ 
{ij}(\lambda),m_{i'j'}(\mu)\}$ and
an R-matrix $R(\lambda,\mu)$ is an operator in $\left(\mathbb{R}^n 
\right )^{\otimes 2}$ depending on parameters $\lambda,\mu$ and solving  
the classical Yang-Baxter equation. We are interested in the bracket  
associated with the {\em trigonometric R-matrix} (for the explicit formula for it, which we will not need, 
see \cite{OPRS}).

Our proof of complete integrability of $T_k$ relies on two facts. One is a well-known statement in the theory of integrable systems: {\em spectral invariants of $M(\lambda)$ are in involution with respect to the Sklyanin bracket}. The second is the result proved in \cite{GSV4}: {\em for any network  on a cylinder, the standard Poisson structure on the space of edge weights induces the trigonometric R-matrix bracket on the space of boundary measurement matrices}. Note that the latter claim, as well as the theorem we are about to formulate, applies not just to acyclic networks on a cylinder but to any network  with sources and sinks belonging to different components of the boundary (in the presence of cycles the definition of boundary measurements has to be adjusted in that path weights may acquire a sign). 

Let $\can$ be a network on the cylinder, and $\bar\can$ be the network  on the torus obtained
via the gluing procedure described in Section~3.1.

\begin{theorem}
\label{cylinder_to_torus} For any network $\hat\can$ on the torus, there exists a network $\can$ on the cylinder with sources and sinks belonging to different components of the boundary such that
$\bar\can$ is equivalent to $\hat\can$, the map $\mathcal{E}_\can \to \mathcal{F}_{\hat{\can}}$ is Poisson with respect to the standard Poisson structures, and spectral invariants of the image $M_\can(\lambda)$ of the boundary measurement map depend only on $\mathcal{F}_{\hat{\can}}$. In particular, 
spectral invariants of  $M_\can(\lambda)$ form an involutive family of functions on $\mathcal{F}_{\hat{\can}}$ with respect to the standard Poisson structure.
\end{theorem}

\begin{figure}[hbtp]
\centering
\includegraphics[height=0.9in]{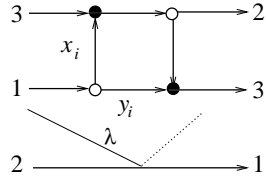}
\caption{Elementary network}
\label{Lax}
\end{figure}

We can now apply the theorem above to networks $\mathcal{N}_{k,n}$. Note that, by construction in Section~3.2, the network $\bar\can_{k,n}$ on the torus is obtained by the gluing procedure from the network  
 $\can_{k,n}$ on the cylinder; the latter is depicted in Fig.~5, provided we refrain from gluing
 together edges marked with the same symbols and regard that figure as representing a cylinder rather than a torus. Furthermore, this network on a cylinder is a concatenation of $n$  {\em elementary networks} of the same form shown on Fig.~\ref{Lax} (for the case $k=3$). 

The boundary measurement matrix that corresponds to the elementary network is
$$
L_i(\lambda)=\left(\begin{array}{cc}
\lambda x_i&x_i+y_i\\
\lambda &1\\
\end{array}\right)
$$
for $k=2$ and
\begin{equation} \label{Laxmat}
L_i(\lambda)=\left(\begin{array}{cccccc}
0&0&0&\dots&x_i&x_i+y_i\\
\lambda&0&0&\dots&0&0\\
0&1&0&\dots&0&0\\
0&0&1&\dots&0&0\\
\dots&\dots&\dots&\dots&\dots&\dots\\
0&0&0&\dots&1&1\\
\end{array}\right)
\end{equation}
for $k\ge 3$. Then the boundary measurement matrix that corresponds to $\mathcal{N}_{k,n}$ is
$$M_{k,n}(\lambda)=L_1(\lambda) \cdots L_n(\lambda)\ .
$$

\begin{proposition} {\rm (i)} There is a Poisson 
bracket on $\mathbb{R}^{2n} = \{ (x_i)_{i=1}^n,  (y_i)_{i=1}^n \}$ 
that induces a trigonometric R-matrix bracket \eqref{sklya} on $M_{k,n}(\lambda)$. For $n\geq 2k-1$, this Poisson bracket is given by
\begin{eqnarray} \nonumber \label{Poiss}
&\{x_i,x_{i+l}\}=-x_ix_{i+l},\ 1\le l \le k-2; \ \{y_i,y_{i+l}\}=-y_iy_{i+l},\ 1\le l\le k-1;\\&
\{y_i,x_{i+l}\}=-y_ix_{i+l},\ 1\le l\le k-1;\ \{y_i,x_{i-l}\}=y_ix_{i-l},\ 0\le l\le k-2\ ,\ \ 
\end{eqnarray}
where indices are understood $\mod n$ and only non-zero brackets are listed.

{\rm (ii)} The bracket \eqref{Poiss} is invariant under the map $T_k$.
\label{bracketxy}
\end{proposition}

\begin{remark}
{\rm {1.} There are formulae similar to  (\ref{Poiss}) for  $T_k$-invariant Poisson bracket in the case
 $n< 2k-1$ as well. Our focus on the ``stable range" $n\geq 2k-1$ will be justified by the geometric interpretation of the maps $T_k$ in Section \ref{geom}. 
 
 {2.} The bracket (\ref{Poiss}) is degenerate.  If $n$ is even and $k$ is odd, there are four Casimir functions :
$
\prod_{i} x_{2i},\quad \prod_{i} x_{2i+1},\quad \prod_{i} y_{2i}, \quad \prod_{i} y_{2i+1}
$. Otherwise, there are two Casimirs : $ \prod x_i$ and $\prod y_i$.
}
\end{remark}

The ring of spectral invariants of $M_{k,n}(\lambda)$ is generated by coefficients of its characteristic polynomial
\begin{equation}
\det(M_{k,n}(\lambda) - z ) = \sum_{i=1}^k \sum_{j=1}^n I_{ij}(x,y) z^i \lambda^j.
\label{integrals}
\end{equation}
(Some of the coefficients $I_{ij}$ are identically zero.)
\begin{theorem}\label{main} Functions $I_{ij}(x,y)$ are preserved by the map $T_k$ and generate
a complete involutive family with respect to the Poisson bracket \eqref{Poiss}. Thus $T_k$ is completely integrable.
\end{theorem}

\begin{corollary} 
Any rational function of $I_{ij}$ that is homogeneous of degree zero
in variables $x,y$, depends only on $p, q$ and is preserved by the map $\overline T_k$. On the level set $\{ \Pi  
p_i q_i =1 \}$, such functions generate
a complete involutive family of integrals for the map $\barT_k$. 
\label{completepq}
\end{corollary}

\begin{remark}
{\rm In general,
these functions define a continuous integrable system on level sets  
of the form $\{ \Pi p_i= c_1, \Pi q_i =c_2 \}$, and the map $\barT_k^c$ intertwines  
the flows of this system on different level sets lying on the same hypersurface $c_1c_2=c$.
Numerical evidence suggests that $\barT_k^c$ is not integrable whenever $c\ne 1$.}
\end{remark}

Involutivity of functions $I_{ij}$ follows from properties of the boundary measurement map mentioned above and Proposition \ref{bracketxy}. The fact that $I_{ij}$ are integrals can be deduced from the fact that the transformation from $M_{k,n}(\lambda)(x,y)$ to $M_{k,n}(\lambda)(T_k(x,y))$ can be described as a composition of a sequence of Postnikov's transformations that do not
affect boundary measurements and a transformation that consists in cutting the cylinder in two and then re-gluing the right boundary of the right  half-cylinder to the left boundary of the left half-cylinder to obtain a new network with the same underlying graph. The latter transformation preserves spectral invariants of the boundary measurement matrix. Another way to see the invariance of $I_{ij}$ under $T_k$, based on a zero-curvature representation, is discussed below.

{\em A zero curvature (Lax) representation\/} with a spectral parameter for a nonlinear dynamical system is a compatibility condition for an over-determined system of linear equations; this is a powerful method of establishing algebraic-geometric complete integrability, see, e.g., \cite{DKN}.  

\begin{proposition}
A zero curvature representation for our discrete-time system is given by
$$
L_i^*(\lambda)= P_i(\lambda) L_{i+r-1}(\lambda) P_{i+1}^{-1}(\lambda)
$$
where the Lax matrix $L_i(\lambda)$ is defined in \eqref{Laxmat}, $L_i^*(\lambda)$ is its image after the transformation $T_k$ 
and $P_i(\lambda)$ is the following auxiliary  matrix
 $$
P_i(\lambda)=\left(\begin{array}{ccccccc}
0&\frac{x_i}{\lambda \sigma_i}&\frac{y_{i+1}}{\lambda \sigma_{i+1}}&0&\dots&0&0\\
0&0&\frac{x_{i+1}}{\sigma_{i+1}}&\frac{y_{i+2}}{\sigma_{i+2}}&\dots&0&0\\
\dots&\dots&\dots&\dots&\dots&\dots&\dots\\
0&0&0&\dots&\frac{x_{i+k-4}}{\sigma_{i+k-4}}&\frac{y_{i+k-3}}{\sigma_{i+k-3}}&0\\
-\frac{1}{\sigma_{i+k-2}}&0&0&\dots&0&\frac{x_{i+k-3}}{\sigma_{i+k-3}}&1\\
\frac{1}{\sigma_{i+k-2}}&-\frac{1}{\lambda \sigma_{i+k-1}}&0&\dots&0&0&0\\
0&\frac{1}{\lambda \sigma_{i+k-1}}&0&\dots&0&0&0\\
\end{array}\right),
$$
where, as before, $\sigma_i=x_i+y_i$. 
\end{proposition}

For $k=3$, one obtains a zero-curvature representation for the pentagram map alternative to the one  given in \cite{So}. The preservation of spectral invariants of $M_{k,n} (\lambda)$
(called, in this context, {\em the monodromy matrix}) follows immediately from the formula
$$
M_{k,n}^*= P_1 L_r \cdots L_n L_1\cdots L_{r-1} P_1^{-1}.
$$

\section{Geometric interpretation} \label{geom}

In this section we give a geometric interpretations of the maps $T_k$. The cases $k\geq 3$ and $k=2$ are treated separately.
\smallskip

{\bf The case $k\geq 3$:  corrugated twisted polygons and higher pentagram maps}. As we already mentioned, a {\it twisted $n$-gon\/} in a projective space is 
 a sequence of points $V_i$ such that $V_{i+n}=M(V_i)$ for all $i \in \Z$ and some fixed projective transformation $M$. The projective group naturally acts on the space of twisted $n$-gons. Let ${\cal P}_{k,n}$ be the 
 space of projective equivalence classes of generic twisted $n$-gons in $\RP^{k-1}$, where ``generic" means that every $k$ consecutive vertices do not lie in a projective subspace. The space ${\cal P}_{k,n}$ has dimension $n(k-1)$. 
 
We say that  a twisted polygon $(V_i)$ is {\it corrugated\/} if, for every $i$, the vertices $V_i,V_{i+1}, V_{i+k-1}$ and $V_{i+k}$ span a projective plane. The projective group preserves the space of corrugated polygons. Denote by ${\cal P}^0_{k,n}\subset \caP_{k,n}$ the space of projective equivalence classes of generic corrugated polygons.
Note that ${\cal P}^0_{3,n}=\caP_{3,n}={\cal P}_n$, the 
space of projective equivalence classes of generic twisted polygons in $\RP^2$.
 The consecutive {\it $(k-1)$-diagonals\/} (the diagonals connecting $V_{i}$ and $V_{i+k-1}$) of a corrugated polygon intersect, and the intersection points form the vertices of a new twisted polygon: the $i$th vertex of this new polygon is the intersection of diagonals  
 $(V_i,V_{i+k-1})$ and $(V_{i+1},V_{i+k})$.
 This $(k-1)$-diagonal map 
 commutes with projective transformations, and hence one obtains a rational map 
 $F_k: {\cal P}^0_{k,n}\to {\cal P}_{k,n}$. 
Note that $F_3$ is the pentagram map; the maps $F_k$ for $k>3$ are called {\it generalized higher pentagram maps}.

We introduce coordinates in a Zariski open subset of ${\cal P}^0_{k,n}$ as follows. 
Our additional genericity assumption is that, for every $i$, 
every three out of the four vertices $V_i, V_{i+1}, V_{i+k-1}$ and $V_{i+k}$ are not collinear.
One can lift the vertices $V_i$ of a corrugated polygon in $\RP^{k-1}$ to vectors $\widetilde V_{i}$  in $\R^k$. Then the four vectors, $\widetilde V_{i}, \widetilde V_{i+1}, \widetilde V_{i+k-1}, \widetilde V_{i+k}$ are linearly dependent for each $i$. 
The lift can be chosen so that the linear recurrence holds
\begin{equation} \label{recurr}
\widetilde V_{i+k}=y_{i-1} \widetilde V_{i} + x_i \widetilde V_{i+1} + \widetilde V_{i+k-1}
\end{equation}
where $x_i$ and $y_i$ are $n$-periodic sequences. Conversely, the recurrence (\ref{recurr}) determines an element in ${\cal P}^0_{k,n}$. 

Recall the notion of the projective duality. Let $P=(V_i)$ be a generic (twisted) polygon in $\RP^{k-1}$. The dual polygon $P^*$ is the polygon in the dual space $(\RP^{k-1})^*$ whose consecutive vertices are the hyperplanes through $(k-1)$-tuples of consecutive vertices of $P$. 

The next theorem describes the geometry of corrugated twisted polygons and interprets the map $T_k$ as 
the generalized higher pentagram  map $F_k$.

\begin{theorem} \label{corrug}
{\rm (i)} The image of a corrugated polygon under the map $F_k$ is a corrugated polygon.

{\rm (ii)} In the $(x,y)$-coordinates, the map $F_k$ is $T_k$, that is, it is
given by the formula  \eqref{mapxy}.

{\rm (iii)} The polygon projectively dual to a corrugated polygon is corrugated.

{\rm (iv)} In the $(x,y)$-coordinates, the projective duality is $(-1)^kD_k$, that is, it is
given, up to a sign, by the formula \eqref{DTD}.
\end{theorem}

\begin{remark}
{\rm A related higher dimensional generalization of the pentagram map was
 considered by M. Glick \cite{Glnot}. An alternative construction of an
 integrable 3-dimensional generalization of the pentagram map is suggested
 by B. Khesin and F. Soloviev \cite{KhSo}.}
\end{remark}

Statement (ii) above, along with Theorem \ref{main}, implies that  the generalized higher pentagram map $F_k$ is completely integrable.

One  also has the $(k-1)$-diagonal map on twisted polygons in the projective plane. We call these maps {\it higher pentagram maps}.
We assume that the polygons are generic in the following sense: for every $i$, the vertices $V_i,V_{i+1},V_{i+k-1}$ and $V_{i+k}$ are in general position, that is, no three are collinear. The $(k-1)$-diagonal map assigns to a twisted $n$-gon $(V_i)$ the twisted $n$-gon whose consecutive vertices are the intersection points of the lines $(V_i,V_{i+k-1})$ and $(V_{i+1},V_{i+k})$. Denote by 
$G_k:{\cal P}_n\to {\cal P}_n$ the respective higher pentagram map. 

Assuming that the above genericity assumption needed for~\eqref{recurr} holds true, one can 
lift the points $V_i$ to vectors $\widetilde V_i\in \R^3$ so that (\ref{recurr}) holds. This defines functions $x_i,y_i$ on a Zariski open subset of ${\cal P}_n$ and provides a map $\psi$ from 
this subset to the $({\bf x,y})$-space. The relation of $G_k$ with $T_k$ is as follows.

\begin{proposition} \label{diags}
{\rm (i)} The map $\psi$ conjugates $G_k$ and $T_k$, that is, $\psi\circ G_k=T_k\circ \psi$.

{\rm (ii)} The map $\psi$ is $k\choose 3$-to-one.
\end{proposition} 

It follows that if $I$ is an integral of the map $T_k$ then $I\circ \psi$ is an integral of the map $G_k$. Thus the higher pentagram maps $G_k$ are integrable. 

\medskip

{\bf The case $k=2$: leapfrog map and circle patterns}. In this case, we are concerned with twisted polygons in $\RP^1$, and we use the following definition of  cross-ratio:
$$
[a,b,c,d]=\frac{(a-b)(c-d)}{(a-d)(b-c)}.
$$

Let ${\cal S}_n$ be the space whose points are pairs of twisted $n$-gons $(S^{-},S)$ in $\RP^1$ with the same monodromy. Here $S$ is a sequence of points $S_i \in \RP^1$, and likewise for $S^{-}$. One has: dim ${\cal S}_n=2n+3$. The group $PGL(2,\R)$ acts on ${\cal S}_n$.  Let $\varphi$ be the map from  ${\cal S}_n$ to the $({\bf x,y})$-space given by the formulas:
\begin{equation*} \label{xyST}
x_i=\frac{(S_{i+1}-S^{-}_{i+2})(S^{-}_i-S^{-}_{i+1})}{(S^{-}_i-S_{i+1})(S^{-}_{i+1}-S^{-}_{i+2})}, \quad y_i=\frac{(S^{-}_{i+1}-S_{i+1})(S^{-}_{i+2}-S_{i+2})(S^{-}_i-S^{-}_{i+1})}{(S^{-}_{i+1}-S_{i+2})(S^{-}_i-S_{i+1})(S^{-}_{i+1}-S^{-}_{i+2})}.
\end{equation*}

\begin{proposition} \label{phi}
{\rm (i)} The fibers of this maps are the $PGL(2,\R)$-orbits, and hence the $({\bf x,y})$-space is identified with the moduli space ${\cal S}_n/PGL(2,\R)$.

 {\rm (ii)} The composition of $\varphi$ with the projection $\pi$ is given by the formulas
$$
p_i=[S^{-}_{i+1},S_{i+1},S^{-}_{i+2},S_{i+2}],\quad  q_i=\frac{[S^{-}_i,S_{i+1},S_{i+2},S^{-}_{i+3}] [S^{-}_{i+1},S^{-}_{i+2},S_{i+2},S^{-}_{i+3}]}{[S^{-}_i,S^{-}_{i+1},S^{-}_{i+2},S^{-}_{i+3}][S^{-}_{i+1},S_{i+1},S_{i+2},S^{-}_{i+3}]}.
$$

 {\rm (iii)} The image of the map $\pi\circ\varphi$ belongs to the hypersurface $\prod p_iq_i=1$.
\end{proposition}

\begin{figure}[hbtp]
\centering
\includegraphics[height=.6in]{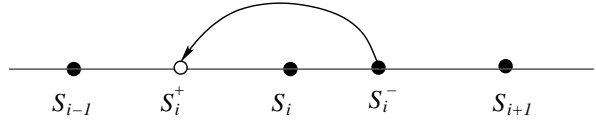}
\caption{Evolution of points in the projective line}
\label{leap}
\end{figure}

Define  a transformation $F_2$ of the space ${\cal S}_n$, acting as $F_2(S^{-},S)=(S,S^{+})$, where $S^{+}$ is given by the following local ``leapfrog" rule: given a quadruple of points $S_{i-1}, S^{-}_i, S_i, S_{i+1}$, the point $S_i^{+}$ is the result of applying to $S^{-}_i$ the unique projective transformation that fixes $S_i$ and interchanges $S_{i-1}$ and $S_{i+1}$. Transformation $F_2$ can be defined this way over any ground field, however in $\RP^1$ we can
interpret  the point $S_i^{+}$ as  the reflection of $S^{-}_i$ in $S_i$ in the projective metric on the segment $[S_{i-1},S_{i+1}]$, see Fig.~\ref{leap}. Recall that the projective metric on a segment is the Riemannian metric whose isometries are the projective transformations preserving the segment, see \cite{Bu}.

\begin{figure}[hbtp]
\centering
\includegraphics[height=2in]{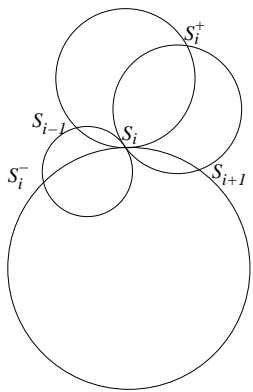}
\caption{Pairwise tangent circles}
\label{circles}
\end{figure}

 Now let the ground field be $\C$, so that the ambient space is $\CP^1$. Define another map, $H_2$, on the space ${\cal S}_n$ by the following local rule. Start with  a quadruple of points $S_{i-1}, S^{-}_i, S_i, S_{i+1}$. Draw the circle through  points $S_{i-1},S^{-}_i,S_{i}$, and then draw the circle through points $S_i, S_{i+1}$, tangent to the previous circle. Now repeat the construction: draw the circle through points $S_{i+1},S^{-}_i,S_{i}$, and then draw the circle through points $S_i, S_{i-1}$, tangent to the previous circle. Finally, define $S^{+}_i$ to be the intersection point of the two ``new" circles, see 
 Fig.~\ref{circles}. This circle pattern generalizes the one studied by O.~Schramm in~\cite{Sc} (there the pairs of circles were orthogonal), see also \cite{BH}.

The next theorem gives geometric interpretations of the map $T_2$. 

\begin{theorem} \label{k2}
{\rm (i)} Over $\C$, the maps $F_2$ and $H_2$ coincide.

 {\rm (ii)} These maps are given by the following equivalent equations:
\begin{equation} \label{Men}
\begin{split}
\frac{1}{S^{+}_i-S_i}+\frac{1}{S^{-}_i-S_i}=\frac{1}{S_{i+1}-S_i}+\frac{1}{S_{i-1}-S_i},\\
\frac{(S^{+}_i-S_{i+1})(S_i-S^{-}_i)(S_i-S_{i-1})}{(S^{+}_i-S_i)(S_{i+1}-S_i)(S^{-}_i-S_{i-1})}=-1.
\end{split}
\end{equation}

{\rm (iii)} The map induced by $F_2$ on the moduli space ${\cal S}_n/PGL(2,\R)$ is the map $T_2$  given in \eqref{mapxy}. 
\end{theorem}

\begin{remark} \label{last}
\rm{ (i) The  closed 2-form on ${\cal S}_n$
$$
\omega=\sum_{i=1}^n \frac{dS^{-}_i\wedge dS_i}{(S^{-}_i-S_i)^2}.
$$
is invariant under the map $F_2$.

(ii) 
Consider the configuration in Fig.~\ref{Menel}. We view the points $S_{i-1}, S^{-}_i,S^{+}_i$ and $S_{i+1}$ as consecutive vertices of a polygon, and the point $S_i$ as the result of the pentagram map (while ${\bar S}_i$ is the result of the inverse pentagram map). By the Menelaus theorem,
$$
\frac{(S^{+}_i-S_{i+1})(S_i-S^{-}_i)({\bar S}_i-S_{i-1})}{(S^{+}_i-{\bar S}_i)(S_{i+1}-S_i)(S^{-}_i-S_{i-1})}=-1.
$$
If the polygonal line $S_{i-1}, S^{-}_i,S^{+}_i,S_{i+1}$ degenerates to a straight line in such a way that the points $S_i$ and ${\bar S}_i$ have the same limiting position then the Menelaus theorem becomes the second equation in (\ref{Men}). The relation of  the Menelaus theorem with  discrete completely integrable systems was discussed in  \cite{KS}.

\begin{figure}[hbtp]
\centering
\includegraphics[height=1.5in]{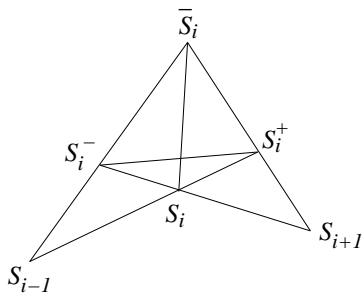}
\caption{Menelaus theorem}
\label{Menel}
\end{figure}

A different construction of a map on twisted polygons in ${\RP}^1$
 that can also be viewed as a limited case of the pentagram map was
 suggested by M. Glick \cite{Glnot}.

(iii) One  interprets equation (\ref{Men}) as a Toda-type equation on the sublattice $\{(m,n)\}$ of $\Z^2$ given by the condition $m+n=0\ ({\rm mod}\ 2)$, see \cite{BH,BS}: 
\begin{equation} \label{our}
\frac{1}{z_{m,n}-z_{m+1,n+1}}+\frac{1}{z_{m,n}-z_{m-1,n-1}}=\frac{1}{z_{m,n}-z_{m+1,n-1}}+\frac{1}{z_{m,n}-z_{m-1,n+1}}.
\end{equation}
Consider another evolution on $\Z^2$, called the cross-ratio equation:
\begin{equation} \label{cross}
[z_{m,n},z_{m+1,n},z_{m+1,n+1},z_{m,n+1}]=q
\end{equation}
where $q$ is a constant. The following two claims hold, see \cite{BH}:\\
a) equation (\ref{cross}), restricted to the sublattice, induces equation (\ref{our});\\
b) given a solution to equation (\ref{our}), a constant $q$,  and a value of $z_{0,1}$, there exists a unique extension to the full lattice $\Z^2$ satisfying (\ref{cross}). The restriction of this extended solution on the other, odd, sublattice also satisfies (\ref{our}).

A continuum limit of equation (\ref{cross}) is known to be the Schwarzian KdV equation \cite{NC}. Thus equation (\ref{our}) also has this continuum limit. Continuum limits of higher  pentagram maps and their bi-Hamiltonian
properties will be discussed in our forthcoming detailed paper.} 
\end{remark}
\bigskip 

{\bf Acknowledgments}. It is a pleasure to thank the Hausdorff Research Institute for Mathematics whose hospitality the authors enjoyed in summer of 2011. We are grateful to A.~Bobenko, V.~Fock, S.~Fomin, M.~Glick, B. Khesin, V.~Ovsienko, R.~Schwartz, F. Soloviev, Yu.~Suris  
for stimulating discussions. M.~G. was partially supported by the
NSF grant DMS-1101462; M.~S. was partially supported by the NSF grants DMS-1101369 and DMS-0800671;
S.~T. was partially supported by the Simons Foundation grant No 209361 and by the NSF grant DMS-1105442;
A.~V. was partially supported by the ISF grant No 1032/08.

\end{document}